# Fourier Series of the Derivatives of Hurwitz and Lerch Zeta Functions

VIVEK  VISHWANATH  RANE

THE INSTITUTE OF SCIENCE,

15,  MADAME  CAMA ROAD,

MUMBAI-400 032 .

INDIA

v_v_rane@yahoo.co.in

**Abstract** :  As a function of second variable , we identify Fourier series of Hurwitz zeta function and its derivatives on unit interval . Consequently , we obtain results based on the formula for Fourier coefficients and also on Parseval's theorem . We do likewise in the case of  Lerch's  zeta function and its derivatives .

# Fourier Series of the Derivatives of Hurwitz and Lerch Zeta Functions

VIVEK VISHWANATH RANE

THE INSTITUTE OF SCIENCE,

15, MADAME CAMA ROAD,

MUMBAI-400 032 , INDIA

v_v_rane@yahoo.co.in

Let $s = \sigma + it$ be a complex variable , where $\sigma$ and t are real . For a complex number $\alpha \neq 0, -1, -2, -3, \ldots\ldots$, let $\zeta(s, \alpha)$ be the Hurwitz zeta function defined by $\zeta(s, \alpha) = \sum_{n \geq 0} (n + \alpha)^{-s}$ for $\sigma > 1$ and its analytic continuation. For an integer $k \geq 0$, let $\zeta_k(s, \alpha) = \sum_{n \geq k} (n + \alpha)^{-s}$ for $\sigma > 1$ and its analytic continuation. Thus $\zeta_k(s, \alpha) = \zeta(s, \alpha + k) = \zeta(s, \alpha) - \sum_{0 \leq n \leq k-1} (n + \alpha)^{-s}$ , where empty sum means zero . It is well-known that $\zeta(s, \alpha)$ ( for a fixed $\alpha$ ) is an analytic function of s except for a pole at $s = 1$. In author [4] , it has been shown that $\zeta(s, \alpha)$ (for a fixed s) is an analytic function of $\alpha$ except for $\alpha = 0, -1, -2, \ldots\ldots$. Note $\zeta(s, \alpha) = \alpha^{-s} + \zeta_1(s, \alpha)$. In author [6] , $\zeta_1(s, \alpha)$ has been shown to be an analytic function of $\alpha$ in the open disc $|\alpha| < 1$. We shall write $\frac{\partial^r}{\partial s^r} \zeta_k(s, \alpha) = \zeta_k^{(r)}(s, \alpha)$ for integers $r, k \geq 0$. In author [3] , we have already defined the Fourier series on an interval of an integrable function . The object of this paper is to identify and study the Fourier series of $\zeta_k^{(r)}(s, \alpha)$ on the interval $[0, 1]$ , as a function of $\alpha$ , for a fixed s in certain right or left half s-plane and also the corresponding results for Lerch's zeta function. This will be done by observing a simple fact that (for a fixed s in a certain region) if the trigonometric or exponential series of functions of $\alpha$ ,

corresponding to the Fourier series on the interval [a,b] of an integrable function (say ,

$\zeta_k^{(r)}(s,\alpha)$) converges uniformly or boundedly on the interval [a,b] , then the said

trigonometric (exponential) series is the Fourier series of that function on the interval .

On similar lines for complex $\alpha \neq 0,-1,-2,......$ , and for real $\lambda$ , we define

Lerch's zeta function $\phi(\lambda,\alpha,s) = \sum_{n \geq 0} e^{2\pi i \lambda n}(n+\alpha)^{-s}$ for $\sigma > 1$ and its analytic continuation. Note that if $\lambda$ is an integer , then $\phi(\lambda,\alpha,s) = \zeta(s,\alpha)$ . If $\lambda$ is not an

integer , then $\phi(\lambda,\alpha,s)$ is an entire function of s . We write

$\frac{\partial^r}{\partial s^r}\phi(\lambda,\alpha,s) = \phi^{(r)}(\lambda,\alpha,s)$ .We shall write I=[0,1] .We shall also study the Fourier series of $\phi(\lambda,\alpha,s)$ as a function of $\alpha$ on the unit interval [0,1] and also as a function of $\lambda$ on

the unit interval I , when s is fixed in a certain right or left half s-plane. From Theorem 2

of author [5] , it is clear that for any integers $k \geq 1$ and $r \geq 0$, the power series in $\alpha$ of

$\zeta_k^{(r)}(s,\alpha)$ can be obtained from the power series in $\alpha$ of $\zeta_k(s,\alpha)$ in the disc $|\alpha| < k$ by

differentiating it $r$ times with respect to s term-by-term so that we have for $|\alpha| < k$ ,

$\zeta_k^{(r)}(s,\alpha) = \frac{\partial^r}{\partial s^r}\zeta_k(s) + \sum_{n \geq 1}\frac{(-\alpha)^n}{n!}\frac{\partial^r}{\partial s^r}\big(s(s+1)..........(s+n-1)\zeta_k(s+n)\big)$ ,

where $\zeta_k(s) = \zeta(s) - \sum_{n \leq k-1}n^{-s}$ so that $\zeta_k^{(r)}(s) = \zeta^{(r)}(s) - \sum_{n \leq k-1}(-1)^r n^{-s}\log^r n$ .

This gives $\zeta^{(r)}(s,\alpha) = \sum_{0 \leq n \leq k-1}(-1)^r(n+\alpha)^{-s}\log^r(n+\alpha)$

$+ \zeta_k^{(r)}(s) + \sum_{n \geq 1}\frac{(-\alpha)^n}{n!}\frac{\partial^r}{\partial s^r}\big(s(s+1)..........(s+n-1)\zeta_k(s+n)\big)$ for $|\alpha| < k$ .

However $k \geq 1$ is arbitrary . This shows that as a function of $\alpha$ , $\zeta^{(r)}(s,\alpha) - \zeta^{(r)}(s)$ is analytic everywhere except possibly at $\alpha = 0,-1,-2,-3,..........$ or equivalently if $s \neq 1$,

$\zeta^{(r)}(s,\alpha)$ as a function of $\alpha$ , is analytic everywhere except possibly at

$\alpha = 0,-1,-2,-3,..........$ Note that $\zeta(s,\alpha) = \alpha^{-s} + \zeta_1(s,\alpha)$ so that

$\zeta^{(r)}(s,\alpha) = (-1)^r \alpha^{-s} \log^r \alpha + \zeta_1^{(r)}(s,\alpha)$ . Also for $s \neq 1$, $\zeta_1(s,\alpha)$ as a function of $\alpha$ ,

is analytic for $|\alpha| < 1$ so that $\zeta_1^{(r)}(s,\alpha)$ as a function of $\alpha$ , is analytic for $|\alpha| < 1$ and in

particular , analytic at $\alpha = 0$ . This means for $s \neq 1$, as a function of $\alpha$ , $\zeta_1^{(r)}(s,\alpha)$ is analytic everywhere except possibly at $\alpha = -1, -2, -3, \ldots \ldots$ In particular , this means as a function of $\alpha$ , $\zeta_1^{(r)}(s,\alpha)$ is analytic on the closed interval $[0,1]$ for $s \neq 1$. By this, we mean $\zeta_1^{(r)}(s,\alpha)$ is analytic in some open set of $\alpha$ -plane containing the interval $[0,1]$, for

$s \neq 1$. Thus for $\operatorname{Re} s < 1, \zeta_1^{(r)}(s,\alpha)$ is an infinitely differentiable function of $\alpha$ on the

interval I . In particular , this means $\zeta_1^{(r)}(s,\alpha)$ is square-integrable in Lebesgue sense on

I for $\sigma < 1$ . Also $\alpha^{-s} \log^r \alpha$ is absolutely Riemann integrable on I as a function of $\alpha$

for $\sigma < 1$ and is square-integrable on I for $\sigma < \frac{1}{2}$ . Thus overall , $\zeta^{(r)}(s,\alpha)$ is an integrable function of $\alpha$ on I for $\sigma < 1$ ; and is square-integrable on I in Lebesgue sense for $\sigma < \frac{1}{2}$ . Thus $\zeta^{(r)}(s,\alpha)$ has trigonometric (exponential) Fourier series as a function of $\alpha$ on I for $\sigma < 1$ ; and satisfies Parseval's theorem on I for $\sigma < \frac{1}{2}$ .

However it is already known that for $\sigma < 1$ and for $0 < \alpha < 1$ ,

$\zeta(s,\alpha) = 2^s \pi^{s-1} \Gamma(1-s) \sum_{n \geq 1} \sin\left(\frac{\pi s}{2} + 2\pi n \alpha\right) n^{s-1}$ , where $\Gamma$ denotes gamma function . Thus

$\zeta(s,\alpha) = 2^s \pi^{s-1} \Gamma(1-s) \sin \frac{\pi s}{2} \cdot \sum_{n \geq 1} n^{s-1} \cos 2\pi n \alpha + 2^s \pi^{s-1} \Gamma(1-s) \cdot \cos \frac{\pi s}{2} \sum_{n \geq 1} n^{s-1} \sin 2\pi n \alpha$

$= \sum_{n \geq 1} \left( b_n(s) \cos 2\pi n \alpha + c_n(s) \sin 2\pi n \alpha \right)$ , say. This is the trigonometric Fourier series of

$\zeta(s,\alpha)$ on I for $\sigma < 1$ . This can also be written as $\zeta(s,\alpha) = \sum_{|n| \geq 1} e^{2\pi i n \alpha} \Gamma(1-s) \cdot (2\pi i n)^{s-1}$ ,

where the principal value of logarithm has been considered . Consequently for any integer $r \geq 0$ and for $\sigma < 1$ , $\zeta^{(r)}(s,\alpha) = \sum\limits_{|n|\geq 1} e^{2\pi i n\alpha} \cdot \frac{\partial^r}{\partial s^r}\Big(\Gamma(1-s)(2\pi i n)^{s-1}\Big)$ .

$$= \sum_{|n|\geq 1} e^{2\pi i n\alpha} \cdot \sum_{\ell=0}^{r} (-1)^{r-\ell}\binom{r}{\ell}(2\pi i n)^{s-1} \cdot \log^{\ell} 2\pi i n \cdot \Gamma^{(r-\ell)}(1-s)$$

$$= \sum_{|n|\geq 1} e^{2\pi i n\alpha} \cdot (2\pi i n)^{s-1} \cdot \left(\sum_{\ell=0}^{r} (-1)^{r-\ell}\binom{r}{\ell}\log^{\ell} 2\pi i n \cdot \Gamma^{(r-\ell)}(1-s)\right).$$

on using Leibnitz's theorem for differentiation . This is the exponential Fourier series of

$\zeta^{(r)}(s,\alpha)$ for Re $s < 1$, as a function of $\alpha$ on I .

In author [2] , we have shown that for $0 < \sigma < 1$ , $\zeta_1(s,\alpha) = \sum\limits_{n=-\infty}^{\infty} a_n(s)e^{2\pi i n\alpha}$ ,

where $a_o(s) = \frac{1}{s-1}$ and $a_n(s) = \int\limits_{1}^{\infty} e^{-2\pi i n u} u^{-s}$ for $n \neq 0$. Writing $x = 1$ and $y = \frac{|t|}{2\pi}$ , there we have shown that $\sum\limits_{n > y+1}\Big|a_n(s)e^{-2\pi i n\alpha}\Big| << x^{-\sigma} \cdot \sum\limits_{n > y+1}\frac{y}{n(n-y)} << x^{-\sigma}\log(y+2)$ and

$\sum\limits_{n=1}^{\infty}\Big|a_n(s)e^{2\pi i n\alpha}\Big| << x^{-\sigma} \cdot \sum\limits_{n=1}^{\infty}\frac{y}{n(n+y)} << x^{-\sigma}$ . Consequently, the series $\sum\limits_{n=-\infty}^{\infty} a_n(s)e^{2\pi i n\alpha}$ of

functions of $\alpha$ on I converges absolutely and uniformly on I and thus is the Fourier series of $\zeta_1(s,\alpha)$ on I .

Next, we state our Propositions . We stress that in the case of $\zeta_k^{(r)}(s,\alpha)$ (for

$k \geq 0$ and $r \geq o$) , we are considering Fourier series as a function of $\alpha$ on I for a fixed $s$

with Re $s < 1$ .

**<u>Proposition</u> 1** : For $0 \leq \alpha \leq 1$ , let $\zeta_1(s,\alpha) = \zeta(s,\alpha) - \alpha^{-s}$ . Then for $0 < \sigma < 1$ , we

have the Fourier series , as a function $\alpha$ on I namely , $\zeta_1(s,\alpha) = \sum\limits_{n=-\infty}^{\infty} a_n(s)e^{2\pi i n\alpha}$ ,

where $a_n(s) = \int\limits_1^\infty e^{-2\pi i n u} u^{-s} du$ for $n \neq 0$ and $a_o(s) = \frac{1}{s-1}$

so that we have $\int\limits_0^1 \zeta_1(s,\alpha) d\alpha = \frac{1}{s-1}$ and $\int\limits_0^1 \zeta_1(s,\alpha) e^{-2\pi i n \alpha} d\alpha = \int\limits_0^1 e^{-2\pi i n u} u^{-s} du$ for $n \neq 0$ .

We also have Parseval's equality , $\int\limits_0^1 \left|\zeta_1(s,\alpha)\right|^2 d\alpha = \sum\limits_{n=-\infty}^\infty \left|a_n(s)\right|^2$ .

**Proposition 2 :** We have for $\sigma < 1$ , the Fourier series for $\zeta(s,\alpha)$ , as a function of $\alpha$

on the unit interval $[0,1]$ , namely $\zeta(s,\alpha) = \Gamma(1-s) \sum\limits_{|n| \geq 1} (2\pi i n)^{s-1} e^{2\pi i n \alpha}$ , so that

$\int\limits_0^1 \zeta(s,\alpha) d\alpha = 0$ and $\int\limits_0^1 \zeta(s,\alpha) e^{-2\pi i n \alpha} d\alpha = (2\pi i n)^{s-1} \cdot \Gamma(1-s)$ for $n \neq 0$

and we have Parseval's equation

$$\int\limits_0^1 \left|\zeta(s,\alpha)\right|^2 d\alpha = \left|\Gamma(1-s)\right|^2 \sum\limits_{|n| \geq 1} \left|(2\pi i n)^{s-1}\right|^2 = 2^{2\sigma-1} \pi^{2\sigma-2} \cosh \pi t \cdot \left|\Gamma(1-s)\right|^2 \cdot \zeta(2-2\sigma)$$

for $\sigma < \frac{1}{2}$ .

Note : In author [3] , we have already defined the concept of 'finite Fourier series at a

rational point of the unit interval I'.

**Corollary** : I) If $\alpha = \frac{a}{q}$ , where $1 \leq a \leq q$ are integers , we have the finite trigonometric

Fourier series , $\zeta(s, \frac{a}{q}) = 2(2\pi q)^{s-1} \Gamma(1-s) \cdot \sum\limits_{r=1}^q \sin\left(\frac{\pi s}{2} + \frac{2\pi r a}{q}\right) \cdot \zeta\left(1-s, \frac{r}{q}\right)$ .

II) We have Parseval's theorem namely , for fixed complex $s_1, s_2$ with Re $s_1$ ,

Re $s_2$ , Re$(s_1 + s_2) < 1$ ,

$\int\limits_0^1 \zeta(s_1,\alpha) \zeta(s_2,\alpha) d\alpha = 2(2\pi)^{s_1+s_2-2} \Gamma(1-s_1) \Gamma(1-s_2) \cdot \cos\frac{\pi}{2}(s_1-s_2) \zeta(2-s_1-s_2)$ .

Note : Mikolas [1] has obtained Parseval's theorem for $\zeta(s,\alpha)$ . However , his approach

is based on Mellin's transformation formula , which is different from ours .

**<u>Proposition</u> 3** : We have for $\sigma < 1$, the Fourier series of $\zeta'(s, \alpha)$ in the form

$$\zeta'(s, \alpha) = \Gamma(1-s) \sum_{|n| \geq 1} e^{2\pi i n \alpha} (2\pi i n)^{s-1} \left( \log 2\pi i n - \frac{\Gamma'}{\Gamma}(1-s) \right)$$

with principal value of logarithm under consideration .

Note : Putting $s = 0$, we get Kummer's result as the following corollary , on noting

$\zeta'(0, \alpha) = \log \frac{\Gamma(\alpha)}{\sqrt{2\pi}}$ . However, Kummer did not identify it as a Fourier series .

**Corollary** : We have the trigonometric Fourier series on $[0,1]$ ,

$$\zeta'(0, \alpha) = \log \frac{\Gamma(\alpha)}{\sqrt{2\pi}} = \sum_{n \geq 1} \frac{\sin 2\pi n \alpha}{\pi n} (\log 2\pi n + \gamma) + \sum_{n \geq 1} \frac{\cos 2\pi n \alpha}{2n}$$

and we have Parseval's equation , $\int_0^1 (\zeta'(0, \alpha))^2 d\alpha = \sum_{n \geq 1} \frac{1}{n^2} \left( \left( \frac{\log 2\pi n + \gamma}{\pi} \right)^2 + \frac{1}{4} \right)$ ,

where $\gamma$ is Euler's constant .

Next , we consider the higher order derivatives of $\zeta(s, \alpha)$ with respect to s . We already know that for $\sigma < 1$ ,

$$\zeta^{(r)}(s, \alpha) = \sum_{|n| \geq 1} e^{2\pi i n \alpha} (2\pi i n)^{s-1} \cdot \sum_{\ell=0}^{r} (-1)^{r-\ell} \binom{r}{\ell} \Gamma^{(r-\ell)}(1-s) \cdot \log^{\ell} 2\pi i n \ .$$

For each $\ell$ with $0 \leq \ell \leq r$ , the series

$(-1)^{r-\ell} \binom{r}{\ell} \Gamma^{(r-\ell)}(1-s) \cdot \sum_{|n| \geq 1} e^{2\pi i n \alpha} (2\pi i n)^{s-1} \log^{\ell} 2\pi i n$ of functions of $\alpha$ on I is boundedly

convergent . Summing up , the series

$\sum_{\ell=0}^{r} (-1)^{r-\ell} \binom{r}{\ell} \Gamma^{(r-\ell)}(1-s) \sum_{|n| \geq 1} e^{2\pi i n \alpha} (2\pi i n)^{s-1} \log^{\ell} 2\pi i n$ of functions of $\alpha$ on I converges

boundedly to $\zeta^{(r)}(s, \alpha)$ and is thus , the Fourier series of $\zeta^{(r)}(s, \alpha)$ for $\sigma < 1$ .

Thus , we get the following Proposition .

**<u>Proposition 4</u>** :  For any integer $r \geq 0$ and for a fixed $s$ with Re $s < 1$ , $\zeta^{(r)}(s,\alpha)$ , as a function of $\alpha$ on the unit interval [0,1] has Fourier series

$$\zeta^{(r)}(s,\alpha) = \sum_{|n| \geq 1} e^{2\pi i n \alpha} (2\pi i n)^{s-1} \cdot \sum_{\ell=0}^{r} (-1)^{r-\ell} \binom{r}{\ell} \Gamma^{(r-\ell)}(1-s) \cdot \log^{\ell} 2\pi i n$$

**<u>Corollary</u>** :  We have for $\sigma < 1$ , $\int_{0}^{1} \zeta^{(r)}(s,\alpha) d\alpha = 0$

and for $n \neq 0$ , we have

$$\int_{0}^{1} \zeta^{(r)}(s,\alpha) e^{-2\pi i n \alpha} d\alpha = (2\pi i n)^{s-1} \sum_{\ell=0}^{r} (-1)^{r-\ell} \binom{r}{\ell} \Gamma^{(r-\ell)}(1-s) \cdot \log^{\ell} 2\pi i n$$

and for $\sigma < \frac{1}{2}$ , we have Parseval's equation

$$\int_{0}^{1} \left| \zeta^{(r)}(s,\alpha) \right|^2 d\alpha$$

$$= e^{-\pi t} \sum_{n \geq 1} (2\pi n)^{2\sigma - 2} \left| \sum_{\ell=0}^{r} (-1)^{r-\ell} \binom{r}{\ell} \Gamma^{(r-\ell)}(1-s) \cdot \left( \log 2\pi n + \frac{\pi i}{2} \right)^{\ell} \right|^2$$

$$+ e^{\pi t} \sum_{n \geq 1} (2\pi n)^{2\sigma - 2} \left| \sum_{\ell=0}^{r} (-1)^{r-\ell} \binom{r}{\ell} \Gamma^{(r-\ell)}(1-s) \cdot \left( \log 2\pi n - \frac{\pi i}{2} \right)^{\ell} \right|^2$$

Exactly on the same lines, we can study Lerch's zeta function $\phi(\lambda,\alpha,s)$ .

For a fixed complex number $\alpha \neq 0, -1, -2, \ldots\ldots\ldots$ and for real $\lambda$ with $0 < \lambda \leq 1$ , and  for a complex number $s = \sigma + it$ , we  have defined $\phi(\lambda,\alpha,s) = \sum_{n \geq 0} e^{2\pi i \lambda n} (n+\alpha)^{-s}$ for $\sigma > 1$ and its analytic continuation . Then $\phi(\lambda,\alpha,s)$ is an entire function of s , unless $\lambda = 1$ in which case $\phi(1,\alpha,s) = \zeta(s,\alpha)$ . Exactly on the same lines as $\zeta(s,\alpha), \phi(\lambda,\alpha,s)$ is an analytic function of $\alpha$ for $\alpha \neq 0, -1, -2, \ldots\ldots\ldots$ as in author [4] . For $\sigma > 0$ and for a fixed $\alpha$ with $0 < \alpha \le 1$, the series $\sum_{n \ge 0} e^{2\pi i \lambda n} (n + \alpha)^{-s}$ of functions of $\lambda$ on the unit interval $[0,1]$ converges boundedly to $\phi(\lambda, \alpha, s)$ and thus $\sum_{n \ge 0} e^{2\pi i \lambda n} (n + \alpha)^{-s}$ is the Fourier series of $\phi(\lambda, \alpha, s)$ as a function of $\lambda$ on the unit interval $[0,1]$. Consequently, we also get Parseval's equation. This we state as follows.

**Proposition 5** : For Re $s > 0$ and for a fixed $\alpha > 0$, Lerch's zeta function $\phi(\lambda, \alpha, s)$ has Fourier series, as a function of $\lambda$ on the unit interval $[0,1]$, namely,

$$\phi(\lambda, \alpha, s) = \sum_{n \ge 0} e^{2\pi i \lambda n} (n + \alpha)^{-s} .$$

**Corollary** : For $\sigma > 0$, we have for $n \ge 0$, $\int_0^1 \phi(\lambda, \alpha, s) e^{-2\pi i \lambda n} d\lambda = (n + \alpha)^{-s}$

and for $n < 0$, we have $\int_0^1 \phi(\lambda, \alpha, s) e^{-2\pi i \lambda n} d\lambda = 0$

and we have Parseval's equation, $\int_0^1 |\phi(\lambda, \alpha, s)|^2 d\lambda = \zeta(2\sigma, \alpha)$ for $\sigma > \frac{1}{2}$.

Next as in the case $\zeta(s, \alpha)$, we have for $\sigma < 1$,

$\phi(\lambda, \alpha, s) = i(2\pi)^{s-1} \Gamma(1-s) \cdot e^{-\frac{\pi i s}{2}} \sum_{n=-\infty}^{\infty} e^{-2\pi i \alpha(n+\lambda)} (n + \lambda)^{s-1}$. However, when $0 \le \sigma < 1$, the equality occurs only for $0 < \alpha < 1$.

Thus $e^{2\pi i \lambda \alpha} \cdot \phi(\lambda, \alpha, s) = i(2\pi)^{s-1} \Gamma(1-s) \cdot e^{-\frac{\pi i s}{2}} \sum_{n=-\infty}^{\infty} e^{-2\pi i n \alpha} (n + \lambda)^{s-1}$ for $\sigma < 1$.

For $\sigma < 1$, both the series $\sum_{n \ge 0} e^{-2\pi i n \alpha} (n + \lambda)^{s-1}$ and $\sum_{n < 0} e^{-2\pi i n \alpha} (n + \lambda)^{s-1}$ of functions of $\alpha$ on I converge boundedly on I for a fixed $\lambda$ with $0 < \lambda < 1$.

Thus $e^{2\pi i \lambda \alpha} \cdot \phi(\lambda, \alpha, s)$, as a function of $\alpha$ on the unit interval $[0, 1]$ has

$i(2\pi)^{s-1} \Gamma(1-s) \cdot \sum_{n=-\infty}^{\infty} e^{-2\pi i n \alpha} \cdot (n + \lambda)^{s-1}$ as the Fourier series when $\sigma < 1$.

This gives the following

**Proposition 6** : For $\sigma < 1$ and for a fixed $\lambda$ with $0 < \lambda < 1$ , as a function of $\alpha$ on the unit interval $[0,1]$ , we have the Fourier series

$$e^{2\pi i \lambda \alpha} \phi(\lambda, \alpha, s) = i(2\pi)^{s-1} \Gamma(1-s) \cdot e^{-\frac{\pi i s}{2}} \sum_{n=-\infty}^{\infty} e^{-2\pi i n \alpha} \cdot (n+\lambda)^{s-1} \quad .$$

**Corollary** : We have for any integer n and for $\sigma < 1$ ,

$$\int_0^1 \phi(\lambda, \alpha, s) e^{2\pi i (n+\lambda)\alpha} d\alpha = i(2\pi)^{s-1} \Gamma(1-s) \cdot (n+\lambda)^{s-1} . e^{-\frac{\pi i s}{2}}$$

and for $\sigma < \frac{1}{2}$ , we have $\int_0^1 |\phi(\lambda, \alpha, s)|^2 d\alpha = (2\pi)^{2\sigma-2} |\Gamma(1-s)|^2 \sum_{n=-\infty}^{\infty} (n+\lambda)^{2\sigma-2}$

$$= (2\pi)^{2\sigma-2} |\Gamma(1-s)|^2 \left( \zeta(2-2\sigma, \lambda) + \zeta(2-2\sigma, 1-\lambda) \right) .$$